\documentclass[reqno,11pt]{article}

\usepackage{amssymb,amsmath}

\usepackage{epsfig}
\usepackage{graphicx,color}

\usepackage[latin1]{inputenc}

\newtheorem{theo}{Theorem}[section]
 \newtheorem{lem}[theo]{Lemma}
\newtheorem{prop}[theo]{Proposition}
\newtheorem{defn}[theo]{Definition}
\newtheorem{rem}[theo]{Remark}

\newtheorem{theorem}{Theorem}[section]

\def\E{\mathbb{E}}

\def\phi{\varphi}

\newcommand {\nn}{\nonumber}

\newcommand {\noi}{\noindent}

\def\mathsf{\bf}
\def\N{\mathbb{N}}

\def\R{\mathbb{R}}
\def\Z{\mathbb{Z}}

\def\E{\mathrm E}
\def\P{\mathrm P}

\def\text{\mbox}

\def\1{{\bf 1}}

\newcommand{\mbf}[1]{\mbox{\boldmath $#1$}}

\newcommand\beqn{\begin{displaymath}}  % no number
\newcommand\eeqn{\end{displaymath}}

\begin{document}

\title{Quasi-MLE for quadratic ARCH model with long memory }
\author{Ieva Grublyt\.e$^{1,2}$,  \ Donatas Surgailis$^2$, \  Andrius \v Skarnulis$^2$   }
\date{\today \\  \small
\vskip.2cm
$^1$ Universit\'e de Cergy Pontoise, D\'epartament de Math\'ematiques, 95302 Cedex, France \\
$^2$ Vilnius University, Institute of Mathematics and Informatics, 08663 Vilnius, Lithuania
}
\maketitle

\begin{abstract}
We discuss parametric quasi-maximum likelihood estimation for quadratic ARCH process
with long memory introduced in Doukhan et al. (2015) and Grublyt\.e and \v Skarnulis (2015)
with  conditional variance %$\sigma^2_t$
given by a strictly positive quadratic form of observable
stationary sequence. %$r_s, s <t$ introduced in Doukhan et al. (2015) and Grublyt\.e and \v Skarnulis (2015).
We prove consistency and asymptotic normality of the corresponding QMLE estimates, including the estimate
of long memory parameter $0< d < 1/2$.
A simulation study of empirical MSE is included.

\end{abstract}

\noi {\it Keywords:} quadratic ARCH process; LARCH model; long memory; parametric estimation;
quasi-maximum likelihood

\section{Introduction}

Recently, Doukhan et al.~\cite{douk2015} and Grublyt\.e and \v Skarnulis~\cite{gru2015} discussed a class of quadratic
ARCH  models of  the form
\begin{eqnarray}\label{genform}
r_t&=&\zeta_t \sigma_t, \qquad
\sigma^2_t\ =\ \omega^2 + \big(a + \sum_{j=1}^\infty b_j  r_{t-j}\big)^2 + \gamma \sigma^2_{t-1},
\end{eqnarray}
where $\gamma, \omega, a, b_j, j\ge 1 $ are real parameters. In \cite{gru2015},
\eqref{genform} was called the Generalized Quadratic ARCH (GQARCH) model.
By iterating the second equation in \eqref{genform},
the squared volatility in \eqref{genform} can be written as a quadratic form
\begin{equation*}\label{sigma2}
\sigma^2_t \ = \ \sum_{\ell =0}^\infty \gamma^\ell \big\{\omega^2 + \big(a + \sum_{j=1}^\infty b_j  r_{t-\ell- j}\big)^2\big\}
\end{equation*}
in lagged variables $r_{t-1}, r_{t-2}, \cdots$, and hence it
represents a particular case of Sentana's~\cite{sen1995}
Quadratic ARCH model with $p=\infty$. The model  (\ref{genform})
includes the classical Asymmetric
GARCH(1,1) process of Engle~\cite{eng1990} and the Linear ARCH (LARCH) model  of Robinson~\cite{rob1991}:
\begin{equation}\label{larch}
r_t\ =\  \zeta_t \sigma_t, \qquad
\sigma_t\ =\  a + \sum_{j=1}^\infty b_j  r_{t-j}.
\end{equation}
Giraitis et al. \cite{gir2000} proved that the squared stationary solution  $\{r^2_t\}$ of the LARCH model in \eqref{larch}
with $b_j$ decaying as $j^{d-1}, 0< d< 1/2 $
may have long memory autocorrelations. For the GQARCH model in (\ref{genform}),
similar results
were established in  \cite{douk2015} and \cite{gru2015}. Namely, assume that
the parameters $\gamma, \omega, a, b_j, j\ge 1 $ in  (\ref{genform}) satisfy
\begin{equation*}\label{Qbc}
b_j \ \sim \ c \, j^{d-1} \qquad (\exists \ 0< d< 1/2, \ c>0)
\end{equation*}
and
\begin{equation*}
\label{4B}
B^2_2 \mu_4 K_4\ <\  1-\gamma, \qquad \gamma \in [0,1), \qquad a \ne 0,
\end{equation*}
where
\begin{equation*}
\mu_4 := \E \zeta^4_0, \qquad B_2 := \sum_{j=1}^\infty b^2_j,
\end{equation*}
and where $K_4$ is the absolute constant from Rosenthal's inequality in \eqref{rosen},  below. Then (see \cite{gru2015}, Thm.5)
there exists a stationary solution of (\ref{genform}) with $\E r^4_t < \infty $ such that
%satisfies
\begin{equation*} \label{cov2r}
{\rm cov}(r^2_0, r^2_t) \  \sim \  \kappa^2_1 t^{2d-1 }, \qquad  t \to \infty
\end{equation*}
%where $\kappa^2_1 :=  \big(\frac{2a \beta}{1-\gamma - B_2} \big)^2 B(d, 1-2d) \E r_0^2 $.
%Moreover,
and
\begin{equation*}\label{lim2r}
n^{-d-1/2} \sum_{t=1}^{[n\tau]} (r^2_t - \E r^2_t) \ \to_{D[0,1]} \  \kappa_2 W_{d +  (1/2)}(\tau),
 \qquad  n \to \infty,
\end{equation*}
where $W_{d +  (1/2)}$ is a fractional Brownian motion with Hurst parameter $H = d + (1/2) \in (1/2,1)$ and
$\kappa_i >0, i=1,2$ are some constants;  $\to_{D[0,1]}$ stands for the weak convergence
in the Skorohod space $D[0,1]$.

As noted in \cite{douk2015}, \cite{gru2015},  the GQARCH model of  (\ref{genform})   and the LARCH model of \eqref{larch}
have similar long memory and leverage properties
and both can be used for modelling of financial data with the above properties. The main disadvantage of the latter model vs. the
former one seems to be the fact that volatility $\sigma_t$ in \eqref{larch}
may assume negative values and is not separated from below by positive constant $c>0$
as in the case of  (\ref{genform}). The standard quasi-maximum likelihood (QML)  approach to estimation of
LARCH parameters is inconsistent and other estimation methods were developed
in Beran and Sch\"utzner~\cite{bera2009}, Francq and Zakoian~\cite{fra2010}, Levine et al.~\cite{lev2009}, Truquet~\cite{tru2014}.

The present paper discusses QML estimation for the 5-parametric
%the
GQARCH model %(\ref{genform})
\begin{equation}\label{sigma3}
\sigma^2_t(\theta) \ = \ \sum_{\ell =0}^\infty \gamma^\ell \big\{\omega^2 + \big(a + c \sum_{j=1}^\infty j^{d-1}  r_{t-\ell- j}\big)^2\big\},
\end{equation}
depending on parameter $\theta = (\gamma, \omega, a, d, c),
0< \gamma <1, \omega >0, a \ne 0, c \ne 0 $ and $d \in (0,1/2)$. The parametric form  $ b_j = c\, j^{d-1} $ of  moving-average coefficients
in \eqref{sigma3}  is the same as in Beran and Sch\"utzner~\cite{bera2009} for the LARCH model. Similarly as in \cite{bera2009} we discuss
the QML estimator $\widehat \theta_n := \arg \min_{\theta \in \Theta} L_n(\theta), \,
L_n(\theta):=\frac{1}{n}\sum_{t=1}^n \Big(\frac{r_t^2}{\sigma_t^2(\theta)} + \log \sigma_t^2(\theta) \Big) $ involving
exact conditional  variance in  \eqref{sigma3}  depending on infinite past $r_s, -\infty < s <t$, and its more realistic version
$\widetilde \theta_n := \arg \min_{\theta \in \Theta} \widetilde L_n(\theta), $ obtained by replacing
the $\sigma^2_t(\theta)$'s in  \eqref{sigma3} by $\widetilde \sigma^2_t(\theta)$ depending only $r_s, 1\le s < t$
(see Sec.~3 for the definition).
It should be noted that the QML function in  \cite{bera2009} is modified to avoid the  degeneracy of $\sigma^{-1}_t$ in \eqref{larch},
by introducing an additional tuning parameter $\epsilon >0$ which affects the performance of the estimator and whose
choice is a non-trivial task. For the GQARCH model \eqref{sigma3} with $\omega >0$
the above degeneracy problem does not occur and we deal with unmodified QMLE in contrast to  \cite{bera2009}. We also note that
our proofs use
different techniques from   \cite{bera2009}. Particularly, the method of  orthogonal Volterra expansions
of the LARCH model used in \cite{bera2009} is not applicable for model \eqref{sigma3}; see (\cite{douk2015}, Example~1).

This paper is organized as follows. Sec.~2 presents some results of \cite{gru2015} about the existence and properties of stationary solution of the
GQARCH equations in \eqref{genform}. In Sec.~3 we define several QMLE estimators of parameter $\theta $ in
\eqref{sigma3}. Sec.~4 presents the main results of the paper devoted to
consistency and asymptotic normality of the QML estimators. Finite sample performance of these estimators is investigated in
the simulation study of Sec.~5. All proofs are relegated to Sec.~6.

\section{Properties of stationary solution}

In this sec. we recall some facts from \cite{gru2015} about stationary solution of \eqref{genform}.
First, we give the definition of it.  Let ${\cal F}_t = \sigma(\zeta_s, s\le t), t \in \Z$ be the sigma-field generated
by $\zeta_s, s\le t$.

\begin{defn} \label{stat} By stationary %$L^2$-
solution
of \eqref{genform}
we mean a stationary and ergodic martingale difference sequence $\{r_t, {\cal F}_t, t \in \Z\}$
with $\E r_t^2 < \infty, \E [r^2_t|{\cal F}_{t-1}] = \sigma^2_t  $
such that  for any $t \in \Z$ the series $X_t := \sum_{s< t} b_{t-s} r_s$
converges in $L^2$, the series $\sigma^2_t = \sum_{\ell=0}^\infty \gamma^\ell (\omega^2 + (a + X_{t-\ell})^2) $
converges in $L^1$ and   \eqref{genform} holds.
\end{defn}

For real $p\ge 2$,  define
\begin{equation}\label{Bpdef}
B_p := \big(\sum_{j=1}^\infty b_j^2\big)^{p/2},
%\begin{cases}
%\sum_{j=1}^\infty |b_j|^p, &0<p < 2, \\
%\big(\sum_{j=1}^\infty b_j^2\big)^{p/2}, &p\ge 2,
%\end{cases}
\qquad
B_{p,\gamma} :=
%\begin{cases}
%B_p/(1-\gamma^{p/2}), &0<p < 2, \\
B_p/(1-\gamma).
%\end{cases}
\end{equation}
We use the following moment inequality.

\begin{prop} \label{Yp} Let $p \ge 2 $ and
$\{Y_j\}$ be a martingale difference sequence such that $\E |Y_j|^p < \infty $;
$\E [Y_j |Y_1, \cdots, Y_{j-1}] = 0, \, j=2,3, \cdots $.
Then
there exists a constant $K_p$ depending only on $p$ and such that
\begin{equation}\label{rosen}
\E \big|\sum_{j=1}^\infty Y_j\big|^p \ \le \ K_p
%\begin{cases}
%\sum_{j=1}^\infty \E |Y_j|^p, &1\le  p \le 2, \\
\big(\sum_{j=1}^\infty (\E |Y_j|^p)^{2/p}\big)^{p/2}.
%\end{cases}
\end{equation}

\end{prop}

\begin{rem} \label{rosenC} {\rm
Inequality \eqref{rosen} is trivial for $p=2, K_2 =1 $. For $p>2$,  \eqref{rosen}
is a consequence of the Burkholder and Rosenthal inequality (see \cite{bur1973}, \cite{ros1970}).
Os\c ekowski \cite{ose2012} proved that $K^{1/p}_p
\le 4 (\frac{p}{4} + 1)^{1/p} \big(1 + \frac{p}{\log (p/2)}\big)$, in particular,
$K_4^{1/4} \le  32.207 $. 
}
\end{rem}

\begin{prop} \label{Xreq} {\rm (\cite{gru2015})}
Let
$\gamma \in [0,1)$ and
$\{\zeta_t\}$ be an i.i.d. sequence with zero  mean and $|\mu|_p :=\E |\zeta_0|^p<\infty$ for some $p\ge 2$.
Assume that
\begin{equation} \label{cQB}
{K_{p} |\mu|_{p} B_{p,\gamma}}   < 1,
\end{equation}
where $B_{p,\gamma}$ is defined in \eqref{Bpdef} and $K_p$ is the absolute constant in \eqref{rosen}.
Then there exists a unique stationary solution $\{r_t\}$ of \eqref{genform} %with $\E |r_t|^p < \infty $
such that the series $X_t = \sum_{j=1}^\infty b_j r_{t-j}$ converges in $L^p$ and
\begin{equation*}\label{X2mom}
\E |r_t|^p \le C(1 + \E |X_t|^{p}) \quad \text{and} \quad
\E |X_t|^{p} \ \le \frac{C B_p }{1 - K_p |\mu|_p B_{p,\gamma}},
\end{equation*}
where $C>0$ is a constant independent of $\{b_j\}, p,  \gamma,$ and the distribution of $\zeta_0$.
Moreover, for $p=2$ condition \eqref{cQB}, or
\begin{equation} \label{cQB2}
B_{2}  = \sum_{j=1}^\infty b^2_j  < 1 - \gamma
\end{equation}
is necessary for the existence of a  stationary $L^2$-solution of \eqref{genform}.

\end{prop}

\section{QML estimators}

The following assumptions on the parametric GQARCH model in \eqref{sigma3} are imposed.

\medskip

\noi  {\bf Assumption (A)} \ $\{\zeta_t\}$ is a standardized i.i.d. sequence with $\E \zeta_t=0,  \E \zeta_t^2=1$.

\medskip

\noi  {\bf Assumption (B)}  \  $\Theta\subset{\R}^5$ is a compact set of parameters $\theta =(\gamma, \omega, a, d, c)$ defined by

\begin{itemize}

\item[(i)] $\gamma \in [\gamma_1, \gamma_2]$ with $0 < \gamma_1 < \gamma_2 <1$;

\item[(ii)] $\omega \in [\omega_1, \omega_2]$ with $0<\omega_1< \omega_2<\infty$;

\item[(iii)] $a\in [a_1, a_2]$ with $-\infty <a_1< a_2< \infty$;

\item[(iv)]  %$b_j=b_j(\beta,d)=\beta j^{d-1}$ where
$d \in [d_1, d_2] $ with $0< d_1 < d_2 < 1/2 $;

\item[(v)]  $c \in [c_1, c_2]$ with $ 0< c_i = c_i(d,\gamma) < \infty, \,  c_1 < c_2 $
such that
%$B_2= (c^2 \sum_{j>0} j^{2(d-1)})$,
$B_2 = c^2 \sum_{j=1}^\infty j^{2(d-1)} <1 -\gamma$ for any
$c \in [c_1,c_2],  \gamma \in [\gamma_1, \gamma_2], d \in [d_1,d_2] $.

\end{itemize}

We assume that the observations $\{r_t,  1\le t \le n\}$ follow
the model in (\ref{genform}) with the true parameter
$\theta_0  =
(\gamma_0, \omega_0, a_0, d_0, c_0)$ belonging to the interior $\Theta_0$ of $\Theta$ in Assumption (B). The restriction
on parameter $c$ in (v) is due to condition  \eqref{cQB2}.
The QML estimator of $\theta \in \Theta $
is defined as
\begin{equation} \label{th*}
\widehat \theta_n := \arg \min_{\theta \in \Theta} L_n(\theta)
\end{equation}
where
\begin{equation} \label{ml}
L_n(\theta)=\frac{1}{n}\sum_{t=1}^n \Big(\frac{r_t^2}{\sigma_t^2(\theta)} + \log \sigma_t^2(\theta) \Big),
\end{equation}
and $\sigma_t^2(\theta) $ is defined in \eqref{sigma3}, viz.,
\begin{eqnarray}\label{Xth}
\sigma_t^2(\theta)&=&\sum_{\ell =0}^\infty \gamma^\ell \big\{\omega^2 + \big(a + c Y_{t-\ell}(d)\big)^2\big\}, \qquad \text{where}  \\
Y_t(d)&:=&\sum_{j=1}^\infty j^{d-1} r_{t-j}.  \nn
\end{eqnarray}
Note the definitions in \eqref{th*}-\eqref{Xth} depend on (unobserved) $r_s, s \le 0$ and
therefore the estimator in \eqref{th*} is usually referred to as the QMLE given infinite past \cite{bera2009}.
 A more realistic version of \eqref{th*} is defined as
\begin{equation} \label{barth}
\widetilde \theta_n := \arg \min_{\theta \in \Theta} \widetilde L_n(\theta),
\end{equation}
where
\begin{eqnarray} \label{barml}
\widetilde L_n(\theta)&:=&\frac{1}{n}\sum_{t=1}^n \Big(\frac{r_t^2}{
\widetilde \sigma_t^2(\theta)} + \log \widetilde \sigma_t^2(\theta) \Big), \quad \text{where} \\
\widetilde \sigma_t^2(\theta)&:=&\sum_{\ell=0}^{t-1} \gamma^\ell \big\{\omega^2 + \big(a + c \widetilde {Y}_{t-\ell}(d)\big)^2\big\}, \quad
\widetilde {Y}_{t}(d)\ :=\ \sum_{j=1}^{t-1} j^{d-1} r_{t-j}. \nn
\end{eqnarray}
Note all quantities in \eqref{barml} depend only on $r_s, 1\le t \le n$, hence \eqref{barth} is called the QMLE given
finite past. The QML functions in \eqref{ml} and \eqref{barml} can be written as
\begin{equation*}
L_{n}(\theta) =\frac{1}{n} \sum_{t=1}^{n}l_{t}(\theta) \quad \text{and} \quad  \widetilde {L}_{n}(\theta) =\frac{1}{n} \sum_{t=1}^{n}\widetilde {l}_{t}(\theta),
\end{equation*}
respectively, where
\begin{equation}\label{ldef}
l_{t}(\theta) := \frac{r_{t}^{2}}{\sigma_{t}^{2}(\theta)}+\log\sigma_{t}^{2}(\theta), \qquad \widetilde {l}_{t}(\theta)
:= \frac{r_{t}^{2}}{\widetilde \sigma_{t}^{2}(\theta)}+\log \widetilde \sigma_{t}^{2}(\theta).
\end{equation}
Finally, following \cite{bera2009} we define a truncated version of \eqref{barth} involving
the last $O(n^\beta)$ quasi-likelihoods  $\widetilde {l}_{t}(\theta), n- [n^\beta] < t \le n $, as follows:
\begin{equation} \label{betath}
\widetilde \theta_n^{(\beta)} := \arg \min_{\theta \in \Theta} \widetilde L_n^{(\beta)}(\theta), \qquad
\widetilde L_n^{(\beta)}(\theta)\  :=\ \frac{1}{[n^\beta]}\sum_{t=n-[n^\beta]+1}^n  \widetilde {l}_{t}(\theta).
\end{equation}
where $0< \beta < 1 $ is a `bandwidth parameter'. Note that for any $t \in \Z$ and $\theta_0  =
( \gamma_0, \omega_0, a_0, d_0, c_0)\in \Theta $,
the random functions $Y_t(d)$ and $\widetilde Y_t(d)$ in \eqref{Xth} and  \eqref{barml}
are infinitely differentiable w.r.t. $d \in (0,1/2)$  a.s. Hence using the explicit form
of $\sigma_t^2(\theta)$ and $\widetilde \sigma_t^2(\theta)$, it follows that
$\sigma_t^2(\theta), \widetilde \sigma_t^2(\theta), {l}_{t}(\theta),  \widetilde {l}_{t}(\theta),
L_{n}(\theta), \widetilde L_n(\theta), \widetilde L_n^{(\beta)}(\theta)$ etc. are all infinitely
differentiable w.r.t. $\theta \in \Theta_0 $ a.s. We use the notation
\begin{equation}\label{Lth}
L(\theta) :=  \E L_n(\theta)  = \E l_{t}(\theta)
\end{equation}
and
\begin{eqnarray}\label{ABmat}
A(\theta)&:=&\E\left[\nabla^T l_{t}(\theta) \nabla l_{t}(\theta)\right] \qquad \text{and} \qquad
B(\theta)\ :=\ \E\left[\nabla^T \nabla l_{t}(\theta)\right],
\end{eqnarray}
where $\nabla  = (\partial/\partial \theta_1, \cdots, \partial/\partial \theta_5)$ and the superscript $T$ stands for
transposed vector. Particularly, $A(\theta)$ and $B(\theta)$ are $5\times 5$-matrices.
By Lemma \ref{lema1},
the expectations in \eqref{ABmat} are well-defined for any $\theta \in \Theta$ under condition $\E r^4_0 < \infty $. We have
\begin{equation}\label{ABeq}
B(\theta) \ = \ \E [\sigma^{-4}_t(\theta) \nabla^T \sigma^2_t (\theta) \nabla \sigma^2_t(\theta)]
\quad \text{and} \quad  A(\theta) = \kappa_4 B(\theta)
\end{equation}
where $\kappa_4 := \E (\zeta^2_0 -1)^2 >0.$

\section{Main results}

Everywhere in this section $\{r_t\}$ is a stationary solution of model \eqref{sigma3} as defined in Definition \ref{stat} and
satisfying Assumptions (A) and (B) of the previous section. As usual, all expectations are taken with respect
to the true value $\theta_0 = (\gamma_0, \omega_0, a_0, d_0, c_0)  \in \Theta_0$, where $\Theta_0  $ is the interior of the parameter set $\Theta \subset \R^5$.

\begin{theorem} \label{thm1}
(i) Let $\E |r_t|^3<\infty$.
Then $\widehat \theta_n$ in \eqref{th*} is a strongly consistent estimator of $\theta_0$, i.e.
$$
\widehat \theta_n \ \stackrel{a.s.}{\to} \  \theta_0. % \quad n \to \infty.
$$

\noi (ii) Let $\E |r_t|^5<\infty$. Then $\widehat \theta_n$ in \eqref{th*} is asymptotically normal:
\begin{equation}\label{clt1}
n^{1/2}\big(\widehat \theta_{n}-\theta_0\big)\ \stackrel{law}{\rightarrow} \ N(0,\Sigma(\theta_0)),  %\quad n \to \infty,
\end{equation}
where $\Sigma(\theta_0) :=B^{-1}(\theta_0) A(\theta_0) B^{-1}(\theta_0) = \kappa_4 B^{-1}(\theta_0) $ and matrices
$A(\theta), B(\theta)$ are defined in \eqref{ABeq}.

\end{theorem}

The following theorem gives asymptotic properties of `finite past' estimators $\widetilde \theta_{n} $ and $\widetilde \theta^{(\beta)}_{n} $
defined in  \eqref{barth} and \eqref{betath}, respectively.

\begin{theorem}\label{thm2}
(i) Let %Assumptions (A) and (B) hold  and
$\E |r_t|^3 < \infty$ and $0<\beta < 1$. Then
$$
\E |\widetilde \theta_{n} - \theta_0| \ \to \ 0 \qquad \text{and} \qquad  \E |\widetilde \theta^{(\beta)}_{n} - \theta_0| \ \to \ 0.
%\quad n \to \infty.
$$

\noi (ii) Let $\E |r_t|^5<\infty$ and $0<\beta<1-2d_0$. Then
\begin{equation}\label{conv}
n^{\beta/2}(\widetilde \theta_n^{(\beta)}-\theta_0) \ \stackrel{law}{\rightarrow} \
N(0, \Sigma(\theta_0)),
\end{equation}
where $\Sigma(\theta_0)$ is the same as in Theorem \ref{thm1}.

%\smallskip

%\noi (iii) {\bf [patikslinti]}

\end{theorem}

The asymptotic results in Theorems \ref{thm1} and \ref{thm2} are similar to the results of (\cite{bera2009}, Thm. 1-4) pertaining
to  the 3-parametric  LARCH model in \eqref{larch} with $b_j =  c j^{d-1}$, except that  \cite{bera2009} deal with a
modified QMLE involving a  `tuning parameter' $\epsilon >0$. Theorems \ref{thm1} and \ref{thm2} are based on
subsequent Lemmas \ref{lema1}-\ref{lema4}
which describe properties of the likelihood processes defined in  \eqref{ml},  \eqref{barml} and \eqref{ldef}. As noted in Sec.~1,
our proofs use different techniques from  \cite{bera2009} which rely on explicit Volterra series representation
of stationary solution of the LARCH model.

\smallskip

For multi-index ${\mbf i} = (i_1, \cdots, i_5) \in \N^5, \, {\mbf i} \ne
{\mbf 0} = (0, \cdots, 0), \,
|{\mbf i}| := i_1 + \cdots + i_5 $,
denote partial derivative
$\partial^{\mbf i} := \partial^{|{\mbf i}|}/ \prod_{j=1}^5 \partial^{i_j} \theta_{i_j} $.

\begin{lem}  \label{lema1}
 Let $\E |r_t|^{2+p} < \infty $, for some integer
$p \ge 1$. Then for any ${\mbf i}\in \N^5, \, 0 < |{\mbf i}| \le p$,
\begin{eqnarray}\label{Esup}
&&\E \sup_{\theta \in \Theta} |\partial^{\mbf i}  l_t(\theta)| < \infty.
\end{eqnarray}
Moreover, if $\E |r_t|^{2+p+\epsilon} < \infty $ for some $\epsilon>0$ and $p\in \N$ then for any ${\mbf i}\in \N^5, \, 0 \le |{\mbf i}| \le p$
\begin{eqnarray}\label{Esupbar}
&&\E \sup_{\theta \in \Theta}|\partial^{\mbf i}  (l_t(\theta) - \widetilde l_t(\theta))| \to 0,  \qquad t \to \infty.
\end{eqnarray}
\end{lem}

\smallskip

\begin{lem} \label{lema2} The function $L(\theta), \theta \in \Theta$ in \eqref{Lth} is bounded and continuous. Moreover, it attains
its unique minimum at $\theta = \theta_0$.
\end{lem}

%\smallskip

\begin{lem} \label{lema3}  Let $\E r^4_0 < \infty $. Then matrices $A(\theta) $ and $B(\theta) $ in \eqref{ABmat}
are well-defined and strictly positive definite for any $\theta \in \Theta$.

\end{lem}

\smallskip

Write $|\cdot |$  for the Euclidean norm in $\R^5 $ and in $\R^5 \otimes \R^5$ (the matrix norm).

\medskip

\begin{lem} \label{lema4}

\noi (i) Let \ $\E |r_t|^3 < \infty $. Then
\begin{eqnarray} \label{convL}
\sup_{\theta \in \Theta}|L_n(\theta)- L(\theta)|\ \stackrel{a.s.} { \to} \ 0   \qquad \text{and} \qquad
\E \sup_{\theta \in \Theta}|L_n(\theta)-\widetilde L_n(\theta)| \ \to \ 0.
\end{eqnarray}
%where $L(\theta)=\E L_n(\theta)=\E l_t(\theta)$;  see \eqref{Lth}.

\medskip

\noi (ii)  Let \ $\E r_t^4 < \infty $.  Then \
$\nabla  L(\theta) = \E \nabla l_t(\theta)$ \ and
\begin{equation} \label{convL2}
\sup_{\theta \in \Theta} |\nabla L_n(\theta)- \nabla L(\theta)| \ \stackrel{a.s.} { \to} \ 0 \qquad \text{and} \qquad
\E \sup_{\theta \in \Theta} |\nabla L_n(\theta)-  \nabla \widetilde L_n(\theta)| \ \to \ 0.
\end{equation}

\smallskip

\noi (iii) Let \ $\E |r_t|^5 < \infty $.  Then \
$\nabla^T \nabla L(\theta) = \E \nabla^T \nabla \ell_t (\theta) = B(\theta) $ \ (see  \eqref{ABmat}) and
\begin{eqnarray} \label{convL3}
&\sup_{\theta \in \Theta} |\nabla^T \nabla L_n(\theta)-  \nabla^T \nabla L(\theta)| \ \stackrel{a.s.} { \to}  \ 0,   \\  %   \qquad \text{and}  \qquad
&\E \sup_{\theta \in \Theta} |\nabla^T \nabla L_n(\theta)- \nabla^T \nabla \widetilde L_n(\theta)| \ \to \ 0. \label{convL4}
\end{eqnarray}

\end{lem}

\section{Simulation study}

In this section we present a short simulation study of the performance of the
QMLE for the GQARCH model in \eqref{sigma3}. % as well as some results from applications of GQARCH model on empirical stock market data.
The GQARCH model in \eqref{sigma3} was simulated with i.i.d. standard normal innovations $\{\zeta_t\}$. The QMLE procedure was evaluated
for medium-term ($n=1000$) and long-term ($n=5000$) samples.
%In order to observe QML estimators' behaviour as sample size $n$ increases, we consider the cases of medium-term ($n=1000$) and long-term ($n=5000$) samples.
The process was generated for $-n \le t \le n$ using the recurrent formula in \eqref{genform} with appropriately truncated sum  $\sum_{j=1}^{\min(n, t+n)} $ and
zero initial condition $\sigma_{-n-1}=0$. The QMLE
estimation used generated time series $r_t, 1\le t \le n$ with $r_t, -n \le t \le 0$ as the pre-sample.
The numerical optimization procedure minimized the QML function:
\begin{equation*}
L_n = \frac{1}{n}\sum_{t=1}^{n}\left(\frac{r_t^2}{\sigma_t^2} + \log\sigma_t^2\right),
\end{equation*}
with
\begin{equation*}
r_t = \zeta_t\sigma_t, \quad \sigma_t^2 = \omega^2 + \big(a + c\sum_{j=1}^{n}j^{d-1}r_{t-j}\big)^2 + \gamma\sigma_{t-1}^2, \quad t = 1,\cdots, n.
\end{equation*}

Finite-sample  performance of  the QML estimator is studied for fixed values of parameters $\gamma_0=0.7, a_0=-0.2, c_0=0.2$  and
different values of
$\omega_0=0.1, 0.01, 0.001$ and the long memory parameter $d_0 = 0.1, 0.2, 0.3, 0.4$. The above choice of $\theta_0 =
(\gamma_0, \omega_0, a_0, c_0, d_0)$ %in our study
can be explained
by the observation that the QML estimation of $\gamma_0, a_0, c_0$ appears to  be more accurate and stable in comparison with
estimation of $\omega_0$ and $d_0$. The very small values of $\omega_0$ in our experiment reflect the fact that in
most real data studied by us, the estimated QML value of $\omega_0$ was less than $0.05$. The presence of $\omega_0 >0$  in
the GQARCH model in \eqref{sigma3}  is very important for consistency of the QMLE procedure, by guaranteeing that $\sigma^2_t(\theta) $ is separated
from zero. A similar role is played by the `tuning parameter'  $\epsilon >0$ in the LARCH estimation in
\cite{bera2009}, except that
%Our impression is that
$\omega_0 >0$
%plays the role of a tuning parameter %  \cite{bera2009} except that it
is estimated in \eqref{sigma3}   and not
{\it ad hoc} imposed as $\epsilon>0$  in \cite{bera2009}.

The numerical QML minimization was performed using the MATLAB language for technical computing, under the
following constraints:
\begin{eqnarray} \label{bnds}
&&0.001 \leq \gamma \le 0.9, \quad 0\le \omega \le 2, \quad -2\le a \le 2, \quad 0\le d\le 0.5, \\
&&(0.05 - \gamma)\vee (\gamma/999) \le c^2 \zeta(2(1-d))  \le (0.99-\gamma)\wedge (99\gamma), \nn
\end{eqnarray}
where $\zeta(z) = \sum_{j=1}^\infty j^{-z} $ is the Riemann zeta function. The last constraint in \eqref{bnds}
guarantees Assumption (B) (v) with
appropriate $0< c_i(d,\gamma), i=1,2$.

The results of the simulation experiment are presented in
Table 1, which shows the sample R(oot)MSEs of the QML
estimates $\widehat \theta_n = ( \widehat \gamma_n, \widehat \omega_n, \widehat a_n, \widehat c_n, \widehat d_n) $ with
$100$ independent replications, for two sample lengths $n = 1000$ and $n=5000$ and
the above choices of $\theta_0 = (\gamma_0, \omega_0, a_0, c_0, d_0). $  Our observations
from Table 1 are summarized below.

%\smallskip

%\noi Conclusions:

\begin{enumerate}

\item   All RMSEs decrease as $n$ increases. The convergence rate of estimates
seems quite good overall. % except maybe for $d_0 = 0.4$ and $\widehat a_n$.

\item  Parameter $\gamma_0$ is estimated rather accurately. E.g., for  $n= 5000$
RMSE($\widehat \gamma_n$) is very stable for all values of $\omega_0$ and $d_0$.

\item The previous conclusion generally applies  also to the QML estimates  $\widehat a_n, \widehat c_n$ and
$\widehat d_n$
except that  their RMSE  markedly increases when $d_0 = 0.4$.

\item  The QML estimate of $\omega_0 \le 0.01$  %is inaccurate and 
seems
to have a `constant' bias   $\approx 0.02\div 0.03$ for all values of  $d_0$ with $n=5000$.% (see Table 2 below).

\end{enumerate}

\section{Proofs}

\noi {\it Proof of Lemma \ref{lema1}.} We use the following (Fa\`a di Bruno)
differentiation rule:
\begin{eqnarray}\label{leib}
\partial^{\mbf i} \sigma^{-2}_t(\theta)&=&\sum_{\nu=1}^{|{\mbf i}|} (-1)^\nu \nu! \,  \sigma^{-2(1+\nu)}_t(\theta)
\sum_{{\mbf j}_1 + \cdots + {\mbf j}_\nu = {\mbf i}} \chi_{{\mbf j}_1, \cdots, {\mbf j}_\nu}
 \prod_{k=1}^\nu  \partial^{{\mbf j}_k} \sigma^{2}_t(\theta), \\
\partial^{\mbf i} \log \sigma^{2}_t(\theta)&=&\sum_{\nu=1}^{|{\mbf i}|} (-1)^{\nu-1} (\nu-1)! \, \sigma^{-2\nu}_t(\theta)
\sum_{{\mbf j}_1 + \cdots + {\mbf j}_\nu = {\mbf i}} \chi_{{\mbf j}_1, \cdots, {\mbf j}_\nu}
 \prod_{k=1}^\nu  \partial^{{\mbf j}_k} \sigma^{2}_t(\theta), \nn
\end{eqnarray}
where the sum $\sum_{{\mbf j}_1 + \cdots + {\mbf j}_\nu = {\mbf i}}$ is taken over decompositions of ${\mbf i}$ into a sum
of $\nu$ multi-indices ${\mbf j}_k \neq {\mbf 0}, k=1, \cdots, \nu$, and $\chi_{{\mbf j}_1, \cdots, {\mbf j}_\nu}  $ is a combinatorial
factor depending only on ${\mbf j}_k, 1\le k \le \nu $.

Let us prove \eqref{Esup}.
We have $|\partial^{\mbf i}  l_t(\theta)| \le r^2_t |\partial^{\mbf i}  \sigma^{-2}_t(\theta)| + |\partial^{\mbf i}  \log \sigma^2_t (\theta)|$.
Hence using \eqref{leib} and the fact that $  \sigma^2_t(\theta) \ge \omega^2/(1-\gamma) \ge \omega_1^2/(1-\gamma_2)> 0$ we obtain
$$
\sup_{\theta \in \Theta}|\partial^{\mbf i}  l_t(\theta)| \le C (r^2_t + 1)\sum_{\nu=1}^{|{\mbf i}| } \sum_{{\mbf j}_1 + \cdots + {\mbf j}_\nu = {\mbf i}}\
\prod_{k=1}^\nu \sup_{\theta \in \Theta} (|\partial^{{\mbf j}_k} \sigma^{2}_t(\theta)|/\sigma_t(\theta)).
$$
Therefore by H\"older's inequality
\begin{eqnarray}\label{EsupH}
\E\sup_{\theta \in \Theta} |\partial^{\mbf i}  l_t(\theta)|&\le&C
(\E (r^2_t +1)^{(2+p)/2})^{2/(2+p)} \nn \\
&\times& \sum_{\nu=1}^{|{\mbf i}|} \sum_{{\mbf j}_1 + \cdots + {\mbf j}_\nu = {\mbf i}}\
\prod_{k=1}^\nu  \E^{1/q_k} \big( \sup_{\theta \in \Theta} |\partial^{{\mbf j}_k} \sigma^{2}_t(\theta)|/\sigma_t(\theta)\big)^{q_k},
\end{eqnarray}
where $\sum_{j=1}^\nu 1/q_j  \le p/(2+p) $.  Note $ |{\mbf i}| = \sum_{k=1}^\nu |{\mbf j}_k| $ and hence the choice $q_k = (2+p)/|{\mbf j}_k| $ satisfies
$\sum_{j=1}^\nu 1/q_j = \sum_{k=1}^\nu |{\mbf j}_k|/(2+p) \le  p/(2+p)$.  Using \eqref{EsupH} and condition $\E |r_t|^{2+p} \le C, $
relation \eqref{Esup} follows from
\begin{eqnarray}\label{Esup1}
\E \sup_{\theta \in \Theta}  \big(|\partial^{{\mbf j}} \sigma^{2}_t(\theta)|/\sigma_t(\theta)\big)^{(2+p)/|{\mbf j}|} &
<&\infty
\end{eqnarray}
for any multi-index ${\mbf j}\in \N^5, \, 1 \le |{\mbf j}| \le p$.

Consider first the case $|{\mbf j}| = 1$, or the partial derivative $\partial_i \sigma_t^2(\theta) = \partial \sigma^2_t(\theta)/\partial \theta_i,
1\le i \le 5$.  We have
\begin{eqnarray}\label{J0}
\partial_i \sigma_t^2(\theta)&=&\begin{cases}
\sum_{\ell =1}^\infty \ell \gamma^{\ell-1} \big\{\omega^2 + \big(a + c Y_{t-\ell}(d)\big)^2\big\},&\theta_i = \gamma,  \\
\sum_{\ell =0}^\infty \gamma^{\ell} 2\omega,&\theta_i = \omega, \\
\sum_{\ell =0}^\infty  \gamma^{\ell} 2\big(a + c Y_{t-\ell}(d)\big),&\theta_i = a,  \\
\sum_{\ell =0}^\infty \gamma^{\ell} 2\big(a + c Y_{t-\ell}(d)\big) Y_{t-\ell}(d), & \theta_i = c, \\
\sum_{\ell =0}^\infty \gamma^{\ell} 2c \big(a + c Y_{t-\ell}(d)\big)\partial_d Y_{t-\ell}(d), & \theta_i = d.
\end{cases}
\end{eqnarray}
We claim that there exist $C >0, 0< \bar \gamma < 1 $ such that
\begin{eqnarray}\label{Juni}
&&\sup_{\theta \in \Theta} \big| \frac{\partial_i \sigma^2_t(\theta)}{\sigma_t(\theta)}\big|
\ \le\ C(1 + J_{t,0} + J_{t,1}), \quad i=1, \cdots, 5, \quad \text{where}  \\
&&J_{t,0}\ := \ \sum_{\ell =0}^\infty \bar \gamma^\ell \sup_{d \in [d_1,d_2]} |Y_{t-\ell}(d)|, \qquad
J_{t,1}\ := \ \sum_{\ell =0}^\infty \bar \gamma^\ell \sup_{d\in [d_1,d_2] } |\partial_d Y_{t-\ell}(d)|. \nn
\end{eqnarray}
Consider \eqref{Juni} for $\theta_i = \gamma$. Using $\ell^2 \gamma^{\ell-2} \le C \bar \gamma^\ell $ for
all $\ell \ge 1, \gamma \in [\gamma_1, \gamma_2] \subset (0,1)  $ and
some $C>0,  0< \bar \gamma <1$ together with Assumption (B) and
 Cauchy inequality, we obtain  $|\partial_\gamma \sigma_t^2(\theta)|/\sigma_t(\theta)
\le \big( \sum_{\ell =1}^\infty \ell^2 \gamma^{\ell-2} \big\{\omega^2 + \big(a + c Y_{t-\ell}(d)\big)^2\big)^{1/2}
\le C (1 + J_{t,0}) $ uniformly in $\theta \in \Theta $, proving  \eqref{Juni} for $\theta_i = \gamma$.
Similarly, $|\partial_c \sigma_t^2(\theta)|/\sigma_t(\theta)
\le C(1+ J_{t,0})$ and $ |\partial_d \sigma_t^2(\theta)|/\sigma_t(\theta)
\le C(1+ J_{t,1})$.  Finally,
for $\theta_i = \omega$ and $\theta_i = a$,  \eqref{Juni} is immediate from \eqref{J0}, proving  \eqref{Juni}.

With  \eqref{Juni} in mind, \eqref{Esup1} for $|{\mbf j}| = 1$  follows from
\begin{eqnarray}\label{Xth3}
\E J^{2+p}_{t,i}\ = \ \E \big(\sum_{\ell =0}^\infty \bar \gamma^\ell \sup_{d \in [d_1,d_2]} |\partial^i_d  Y_{t-\ell}(d)|\big)^{2+p} \ < \ \infty, \quad i=0,1.
\end{eqnarray}
Using Minkowski's inequality and stationarity of $\{Y_t(d)\}$ we obtain $\E^{1/(2+p)} J^{2+p}_{t,i}
\le \sum_{\ell =0}^\infty \bar \gamma^{\ell} $  $ \E^{1/(2+p)} \sup_{d} |\partial^i_d Y_{t-\ell}(d)|^{2+p}  \le
C ( \E\sup_{d} |\partial^i_d Y_{t}(d)|^{2+p})^{1/(2+p)} $, where $\partial^i_d Y_{t}(d) = \sum_{j=1}^\infty \partial^i_d j^{d-1} r_{t-j}. $
Hence using (\cite{bera2009}, Lemma 1 (b)) and the inequality $x y \le x^q/q + y^{q'}/q', x,y >0, 1/q + 1/q' =1 $
we obtain
\begin{eqnarray}
\sum_{i=0}^1\E J^{2+p}_{t,i}
&\le&C\sum_{i=0}^1\E\sup_{d\in [d_1,d_2]} |\partial^i_d Y_{t}(d)|^{2+p}\nn \\
&\le&C\sum_{i=0}^2 \sup_{d\in [d_1,d_2]} \E |\partial^i_d Y_{t}(d)|^{2+p} \ < \ \infty \label{Jbdd}
\end{eqnarray}
since $\sup_{d \in [d_1,d_2]}\E |\partial^i_d Y_{t}(d)|^{2+p} \le C \sup_{d \in [d_1,d_2]}
\big(\sum_{j=1}^\infty (\partial^i_d j^{d-1})^2 (\E |r_{t-j}|^{2+p})^{2/(2+p)} \big)^{(2+p)/2}
< \infty $ according to condition $\E |r_t|^{2+p} < C, $ Rosenthal's inequality in \eqref{rosen} and the fact that
$\sup_{d \in [d_1,d_2]} \sum_{j=1}^\infty (\partial^i_d j^{d-1})^2 $  $\le \sup_{d \in [d_1,d_2]} \sum_{j=1}^\infty j^{2(d-1)} (1 + \log^2 j)^2  < C, i=0,1,2$.  This proves \eqref{Esup1} for $|{\mbf j}| = 1 $.

The proof of  \eqref{Esup1} for $2 \le |{\mbf j}| \le p $ is simpler since it reduces to
\begin{eqnarray}\label{Esup2}
\E \sup_{\theta \in \Theta} |\partial^{{\mbf j}} \sigma^{2}_t(\theta)|^{(p+2)/2} &<&\infty, \qquad 2 \le |{\mbf j}| \le p.
\end{eqnarray}
Recall $\theta_1 = \gamma$ and  ${\mbf j}' :=
{\mbf j} - (j_1, 0,0,0,0) = (0,j_2,j_3,j_4,j_5)$. If ${\mbf j}'
= {\mbf 0}$ then $ \sup_{\theta \in \Theta} |\partial^{{\mbf j}} \sigma^{2}_t(\theta)| \le CJ_{t,0}$ follows as in \eqref{Juni} implying
 \eqref{Esup2} as in \eqref{Jbdd} above. Next, let  ${\mbf j}'
\neq {\mbf 0}$. Denote
\begin{equation}\label{Qdef}
Q_t^2(\theta):=  \omega^2 + \big(a + c Y_t(d)\big)^2
\end{equation}
so that $\sigma^2_t(\theta) = \sum_{\ell=0}^\infty \gamma^\ell Q^2_{t-\ell}(\theta)$.
We have with $m: =j_1 \ge 0$ that
$|\partial^{\mbf j} \sigma_t^2(\theta)|
\le\sum_{\ell =m}^\infty (\ell !/(\ell-m)!)\gamma^{\ell-m}   $   $| \partial^{\mbf j'} Q_{t-\ell}^2(\theta)|$ and
% & \le &\sum_{\ell =m}^\infty \frac{\ell !}{m! (\ell-m)!} \gamma^{\ell/2-m} |\partial^{\mbf j'} Q_{t-\ell}^2(\theta)/Q_{t-\ell}(\theta)|
\eqref{Esup1} follows from
\begin{eqnarray}\label{Esup3}
\E \sup_{\theta \in \Theta} |\partial^{{\mbf j}} Q^{2}_t(\theta)|^{(p+2)/2}&<&\infty.
\end{eqnarray}
For $j_2 \ne 0$ (recall $\theta_2 = \omega$) the derivative in \eqref{Esup3}  is trivial so that it suffices to check
\eqref{Esup3} for $j_1 = 0$ only. Then applying Fa\`a di Bruno's rule
we get
$$
|\partial^{{\mbf j}} Q^{2}_t(\theta)|^{(p+2)/2} \le  C\sum_{{\mbf j}_1 + {\mbf j}_2 = {\mbf j}} |\partial^{{\mbf j}_1} (a+ c Y_t(d))|^{(p+2)/2}
|\partial^{{\mbf j}_2} (a+ c Y_t(d))|^{(p+2)/2}
$$
and hence \eqref{Esup3} reduces to
\begin{eqnarray*}\label{Esup4}
\E \sup_{\theta \in \Theta}  |\partial^{{\mbf j}} (a+ c Y_t(d))|^{p+2}&<&\infty, \qquad 0 \le |{\mbf j}| \le p,
\end{eqnarray*}
whose proof is similar to \eqref{Xth3} above. This ends the proof of \eqref{Esup}.

\smallskip

The proof of \eqref{Esupbar} is similar.
We have $|\partial^{\mbf i}  (l_t(\theta) - \widetilde l_t(\theta))|
\le r^2_t |\partial^{\mbf i}(\sigma^{-2}_t(\theta)  -  \widetilde \sigma^{-2}_t(\theta)) | +
|\partial^{\mbf i} (\log \sigma^2_t (\theta)  -  \log \widetilde \sigma^{2}_t(\theta))|$.
%Note the bounds in \eqref{EsupH} and \eqref{Esup1} hold with $l_t(\theta), \sigma^2_t (\theta)$ replaced by
%$\bar l_t(\theta), \bar \sigma^2_t (\theta)$, respectively. For bounded r.v. $\xi$, let $\|\xi\|_\infty := \lim_{p\to \infty} \E^{1/p} |\xi|^p$.
Hence, using H\"older's inequality similarly as in the proof \eqref{Esup}  it suffices to show
\begin{eqnarray} \label{diffbdd1}
\E \sup_{\theta \in \Theta}  |\partial^{\mbf i}(\sigma^{-2}_t(\theta)  -  \widetilde \sigma^{-2}_t(\theta)) | ^{\frac{p+2}{p}}\to 0 \quad \text{and} \quad
\E \sup_{\theta \in \Theta}  |\partial^{\mbf i}(\log \sigma^2_t (\theta)  - \log \widetilde \sigma^{2}_t(\theta))|^{\frac{p+2}{p}}\to 0.
\end{eqnarray}
Below, we prove the first relation in \eqref{diffbdd1} only, the proof of the second one being analogous.

Using the differentiation rule in \eqref{leib} we have that
\begin{eqnarray*}
|\partial^{\mbf i} (\sigma^{-2}_t(\theta)-\widetilde  \sigma^{-2}_t(\theta))|&\le&C \sum_{\nu=1}^{|{\mbf i}|}
\sum_{{\mbf j}_1 + \cdots + {\mbf j}_\nu = {\mbf i}} \big|W_t^{{\mbf j}_1, \cdots, {\mbf j}_\nu}
(\theta) - \widetilde W_t^{{\mbf j}_1, \cdots, {\mbf j}_\nu}(\theta)  \big|,
\end{eqnarray*}
where
\begin{eqnarray*}
&W_t^{{\mbf j}_1, \cdots, {\mbf j}_\nu}(\theta)
:=\sigma^{-2(1+\nu)}_t(\theta)\prod_{k=1}^\nu \partial^{{\mbf j}_k} \sigma^{2}_t(\theta), \\
&\widetilde W_t^{{\mbf j}_1, \cdots, {\mbf j}_\nu}(\theta)
:= \widetilde \sigma^{-2(1+\nu)}_t(\theta)\prod_{k=1}^\nu \partial^{{\mbf j}_k} \widetilde \sigma^{2}_t(\theta).
\end{eqnarray*}
Whence, \eqref{diffbdd1} follows from
\begin{equation}
\sup_{\theta \in \Theta} |W_t^{{\mbf j}_1, \cdots, {\mbf j}_\nu}
(\theta) - \widetilde W_t^{{\mbf j}_1, \cdots, {\mbf j}_\nu}(\theta)|\ \to_p \  0, \qquad t \to \infty \label{V1}
\end{equation}
and
\begin{equation}
\E \sup_{\theta \in \Theta}  \big(|W_t^{{\mbf j}_1, \cdots, {\mbf j}_\nu}
(\theta)| + |\widetilde W_t^{{\mbf j}_1, \cdots, {\mbf j}_\nu}(\theta)|\big)^{(p+2+\epsilon)/p}\ \le \  C < \infty \label{V2}
\end{equation}
for some constants $\epsilon>0$ and $C > 0$ independent of $t$. In turn, \eqref{V1} and \eqref{V2} follow from
\begin{equation}\label{V11}
\sup_{\theta \in \Theta} |\partial^{{\mbf j}} (\sigma^{2}_t(\theta) - \widetilde \sigma^2_t(\theta))| \ \to_{\rm p} \ 0,  %\qquad  \forall \,  |{\mbf j}| \ge 0,
 \qquad t \to \infty
\end{equation}
and
\begin{eqnarray}\label{V21}
\E \sup_{\theta \in \Theta}  \big(|\partial^{{\mbf j}} \sigma^{2}_t(\theta)|/\sigma_t(\theta)\big)^{(2+p+ \epsilon)/|{\mbf j}|} &
<&C, \\
\E \sup_{\theta \in \Theta}  \big(|\partial^{{\mbf j}} \widetilde \sigma^{2}_t(\theta)|/\widetilde \sigma_t(\theta)\big)^{(2+p+ \epsilon)/|{\mbf j}|} &
<&C, \nn
\end{eqnarray}
for any multi-index
${\mbf j}$ such that $  |{\mbf j}| \ge 0 $ and $1 \le |{\mbf j}| \le p$, respectively.

Using condition $\E |r_t|^{2+ p+\epsilon}  < C $,
relations in \eqref{V21}  can be proved analogously to \eqref{Esup1} and we omit the details.
Consider \eqref{V11}. Split $\sigma_t^2(\theta) -  \widetilde \sigma^2_t(\theta) = U_{t,1}(\theta)  + U_{t,2}(\theta)$, where
\begin{eqnarray} \label{Udef}
U_{t,1}(\theta)&:=&\sum_{\ell =1}^{t-1}  \gamma^\ell \big\{\big(a + c Y_{t-\ell}(d)\big)^2  - \big(a + c \widetilde Y_{t-\ell}(d)\big)^2  \big\}, \\
U_{t,2}(\theta)&:=&\sum_{\ell =t}^{\infty}  \gamma^\ell \big\{\omega^2 + \big(a + c  Y_{t-\ell}(d)\big)^2\big\}. \nn
\end{eqnarray}
Then $\sup_{\theta \in \Theta} |\partial^{{\mbf j}}U_{t,i}(\theta)| \to_{\rm p} \ 0,  t\to \infty,  i=1,2 $ follows by using Assumption (B)
and considering the bounds on the derivatives as in the proof of \eqref{Esup1}. For instance, let us prove \eqref{V11} for
$\partial^{{\mbf j}} = \partial_d,
|{\mbf j}| =1 $.  We have
$|\partial_d U_{t,1}(\theta)|  \le   C \sum_{\ell =1}^{t-1} \gamma^\ell \big\{
(1 + |\bar Y_{t-\ell}(d)|) |\partial_d (Y_{t-\ell} (d) -  \widetilde Y_{t-\ell}(d))| +
|\partial_d Y_{t-\ell} (d)|\,|Y_{t-\ell} (d) -  \widetilde  Y_{t-\ell}(d)|\big\}. $
Hence, $\sup_{\theta \in \Theta} |\partial_d U_{t,1}(\theta)| $  $\to_{\rm p} 0 $ follows from
$0\le \gamma \le \gamma_2 < 1 $ and
\begin{eqnarray}
&&\E \sup_{d\in [d_1,d_2]} (|Y_t (d) - \widetilde Y_{t}(d)|^2 +  |\partial_d (Y_{t} (d) -  \widetilde Y_{t}(d))|^2) \ \to \  0  \quad \text{and} \label{Y1} \\
&&\E \sup_{d\in [d_1,d_2]} (|Y_t(d)|^2 + |\widetilde Y_t(d)|^2 +  |\partial_d Y_t (d)|^2 +  |\partial_d \widetilde Y_t(d)|^2) \ \le \  C.  \label{Y2}
\end{eqnarray}
The proof of  \eqref{Y2} mimics that of \eqref{Jbdd} and therefore is omitted. To show \eqref{Y1}, note
$Y_t(d) - \widetilde Y_t(d) = \sum_{j=t}^\infty j^{d-1} r_{t-j} $ and use a similar  argument as
in \eqref{Jbdd} to show that  the l.h.s. of \eqref{Y2} does not exceed
$C \sup_{d \in [d_1,d_2]} \sum_{i=0}^2 \E |\partial^i_d (Y_t(d)-\widetilde Y_t(d))|^2
\le C \sup_{d \in [d_1,d_2]} \sum_{j=t}^\infty j^{2(d-1)} $ $ (1 + \log^2 j) \to 0 \, (t\to \infty)$.
This proves \eqref{V11} for $ |{\mbf j}| =1$ and $\partial^{{\mbf j}} = \partial_d $.
The remaining cases in \eqref{V11}
follow similarly and we omit the details.  This proves \eqref{Esupbar} and completes the proof of
Lemma \ref{lema1}. \hfill  $\Box$

\medskip

\noi {\it Proof of Lemma \ref{lema2}.} We have $|L(\theta_1) - L(\theta_2)|  \le \E |l_{t}(\theta_1) -   l_{t}(\theta_2)|
\le C \E |\sigma^2_t(\theta_1) - \sigma^2_t(\theta_2)|,  $  where the last expectation
can be easily shown to vanish as $|\theta_1 - \theta_2| \to 0, \, \theta_1, \theta_2 \in \Theta$. This proves the first statement
of the lemma. To show the second statement of the lemma,  write
$$
L(\theta) - L(\theta_0) =  \E \big[ \frac{\sigma^2_t(\theta_0)}{\sigma^2_t(\theta)}- \log \frac{\sigma^2_t(\theta_0)}{\sigma^2_t(\theta)} -1 \big].
$$
The function $f(x) := x - 1 - \log x >0$ for $x >0, x \ne 1 $ and $f(x) = 0$ if and only if $x=1$. Therefore
$L(\theta) \ge L(\theta_0), \forall \, \theta \in \Theta $ while $L(\theta) = L(\theta_0) $ is equivalent to
\begin{equation}\label{ss}
\sigma^2_t(\theta) = \sigma^2_t(\theta_0)  \qquad (\P_{\theta_0}-\text{a.s.})
\end{equation}
Thus, it remains to show that \eqref{ss} implies $\theta = \theta_0 = (\gamma_0, \omega_0,a_0,d_0, c_0)$.
Consider the `projection' %{\bf [cia mes nereikalaujame kad  dydziai \eqref{ss} turi antra momenta]}
$P_s \xi = \E [\xi|{\cal F}_s] -  \E [\xi|{\cal F}_{s-1}] $ of r.v. $\xi, \E |\xi| < \infty $, where
$ {\cal F}_s = \sigma(\zeta_u, u \le s)$ (see sec.2). \eqref{ss} implies
\begin{equation}\label{ssp}
0 = P_s(\sigma^2_t(\theta) - \sigma^2_t(\theta_0)) = P_s  (Q_t^2(\theta)-Q_t^2(\theta_0)) + (\gamma - \gamma_0) P_s  \sigma_{t-1}^2(\theta_0),
\qquad \forall \  s \le t-1,
\end{equation}
where $Q_t^2(\theta) = \omega^2 + \big(a + \sum_{u<t} b_{t-u}(\theta) r_u\big)^2$ is the same as in \eqref{Qdef}.
We have
\begin{eqnarray}\label{PQ}
P_sQ_t^2(\theta)
&=&2ab_{t-s}(\theta)r_s + 2b_{t-s}(\theta)r_s \sum_{u<s} b_{t-u}(\theta) r_u + \sum_{s \le u <t} b^2_{t-u}(\theta) P_s r^2_u  \\
&=&2ab_{t-s}(\theta)\zeta_s \sigma_s(\theta_0) + 2b_{t-s}(\theta)\zeta_s \sigma_s(\theta_0)
\sum_{u<s} b_{t-u}(\theta) r_u \nn \\
&+&\sum_{s < u <t} b^2_{t-u}(\theta) P_s \sigma^2_u(\theta_0) + b^2_{t-s}(\theta) (\zeta^2_s-1) \sigma^2_s(\theta_0). \nn
\end{eqnarray}
Whence and from \eqref{ssp} for $s=t-1$ using $P_{t-1}  \sigma_{t-1}^2(\theta_0) = 0$ we obtain
\begin{equation}\label{AB}
C_1(\theta,\theta_0) \zeta^2_{t-1} + 2C_2(\theta,\theta_0) \zeta_{t-1} - C_1(\theta,\theta_0) = 0
\end{equation}
where
\begin{eqnarray*}
C_1(\theta,\theta_0)&:=&(c^2 - c^2_0) \sigma_{t-1}(\theta_0), \\
C_2(\theta,\theta_0)&:=&(ac - a_0 c_0) + \sum_{u< t-1} (c^2 (t-u)^{d-1} - c_0^2 (t-u)^{d_0-1}) r_u.
\end{eqnarray*}
Since $C_i(\theta,\theta_0), i=1,2$ are ${\cal  F}_{t-2}$-measurable,  \eqref{AB} implies
$C_1(\theta,\theta_0)= C_2(\theta,\theta_0) = 0$, particularly, $c = c_0$ since $ \sigma_{t-1}(\theta_0) \ge \omega >0$. Then
$0 =   C_2(\theta,\theta_0) = c_0(a-a_0) + c_0^2  \sum_{u< t-1} ((t-u)^{d-1} - (t-u)^{d_0-1}) r_u $ and $\E r_u = 0$ lead
to $a= a_0$ and next to $0 = \E (\sum_{u< t-1} ((t-u)^{d-1} - (t-u)^{d_0-1}) r_u )^2 =
\E r^2_0 \sum_{j\ge 2} (j^{d-1} - j^{d_0-1})^2 = 0$, or $d = d_0$. Consequently,
$ P_s  (Q_t^2(\theta)-Q_t^2(\theta_0))  = 0$ for  any $s \le t-1 $ and hence
$\gamma = \gamma_0$ in view of \eqref{ssp}.  Finally,  $\omega=\omega_0$ follows from $\E \sigma_{t}^{2}(\theta) = \E\sigma_{t}^{2}(\theta_0)$
and the fact that  $\omega >0, \omega_0>0$.  This proves
$\theta = \theta_0$ and the lemma, too. \hfill $\Box$

\medskip

\noi {\it Proof of Lema \ref{lema3}.} From \eqref{ABeq}, it suffices to show that
\begin{equation} \label{nabS}
\nabla \sigma_t^2 (\theta)\lambda^T  =0
\end{equation}
for some $\theta \in \Theta$ and $\lambda \in \R^5, \lambda \ne 0$ leads to a contradiction. To the last end,
we use a similar projection argument as in the proof of Lemma \ref{lema2}. First, note that
$\sigma_t^2 (\theta) = Q_t^2 (\theta) + \gamma \sigma_{t-1}^2 (\theta) $ implies
$$
\nabla \sigma_t^2(\theta) =(0, \nabla_4 Q_t^2(\theta))+
\gamma \nabla \sigma_{t-1}^2(\theta)+
( \nabla \gamma) \sigma_{t-1}^2(\theta),
$$
where $\nabla_4 = (\partial/\theta_2, \cdots, \partial \theta_5) $. Hence and using the fact that \eqref{nabS} holds for any $t \in \Z$ by stationarity,
from \eqref{nabS} we obtain
\begin{equation} \label{NQ}
(\sigma_{t-1}^2(\theta), \nabla_4 Q_t^2(\theta))\lambda^T  =  0.
\end{equation}
%Since $\lambda^T \dot \sigma_t^2=0$  all $t$ the equation above becomes
%$0=((\dot Q_t^2)^T, 0)\lambda+(0,0,0,0, \sigma_{t-1}^2)\lambda.$
%Moreover, by taking projections of both sides we obtain
Thus,
\begin{equation*}\label{proeq1}
(P_s \sigma_{t-1}^2(\theta), P_s \nabla^T_4 Q_t^2(\theta))\lambda \ =\  0,  \qquad \forall \, s\le t-1;
\end{equation*}
c.f. \eqref{ssp}. For $s = t-1$
using $P_{t-1} \sigma_{t-1}^2(\theta) = 0, \,
P_{t-1} \nabla_4 Q_t^2(\theta) = \nabla_4 P_{t-1} Q_t^2(\theta) $ by differentiating \eqref{PQ}
similarly to \eqref{AB} we obtain
\begin{equation}\label{AD}
D_1(\lambda) \zeta^2_{t-1} + 2D_2(\lambda) \zeta_{t-1} - D_1(\lambda) = 0
\end{equation}
where $D_1(\lambda):= 2\lambda_5\sigma_{t-1}(\theta)$ and
\begin{eqnarray*}
&D_2(\lambda)\ :=\ \lambda_3 c  + \lambda_5 a + 2 \lambda_5 c   \sum_{u< t-1}  (t-u)^{d-1} r_u
+ \lambda_4 c^2  \sum_{u< t-1}  (t-u)^{d-2} \log (t-u) r_u,
\end{eqnarray*}
$\lambda  = (\lambda_1, \cdots, \lambda_5)^T$. As in  \eqref{AB}, $D_i(\lambda), i=1,2 $ are ${\cal F}_{t-2}$-measurable,
\eqref{AD} implying $D_i(\lambda) = 0, i=1,2$. Hence, $\lambda_5 = 0$ and then  $ D_2(\lambda) = 0 $ reduces
to $\lambda_3 c  + \lambda_4 c^2  \sum_{u< t-1}  (t-u)^{d-2} \log (t-u) r_u = 0$. By taking expectation
and using $c \ne 0$ we get $\lambda_3 =  0$ and then $\lambda_4 = 0$ since $\E ( \sum_{u< t-1}  (t-u)^{d-2} \log (t-u) r_u)^2 \ne 0$.
The above facts allow to rewrite  \eqref{NQ} as $2\omega \lambda_2 + \lambda_1 \sigma^2_{t-1}(\theta) = 0$. Unless both $\lambda_1, \lambda_2$ vanish,
the last equation means that either
$\lambda_1 \ne 0$ and
$\{\sigma^2_t (\theta)\}$ is a deterministic process which contradicts $c \ne 0$, or $\lambda_1 = 0, \lambda_2 \ne 0$ and $\omega=0$, which contradicts
$\omega \ne 0$. Lemma \ref{lema3} is proved. \hfill  $\Box$

\medskip

\noi{\it Proof of Lemma \ref{lema4}.} Consider the first relation in \eqref{convL}. The pointwise convergence $L_n(\theta) \stackrel{a.s.}
 { \to} L(\theta) $  follows by ergodicity of $\{l_t(\theta)\}$ and the uniform
convergence in \eqref{convL}
from $\E \sup_{\theta \in \Theta} |\nabla l_t(\theta)| < \infty $, c.f. (\cite{bera2009},  proof of Lemma 3),
which in turn follows from of Lemma \ref{lema1} \eqref{Esup}
with $p=1$.
The proof of the second relation in \eqref{convL} is immediate from   Lemma \ref{lema1} \eqref{Esupbar} with $p=0, \epsilon = 1$.
The proof of the statements
(ii) and (iii) using Lemma \ref{lema1} is similar and is omitted. \hfill $\Box$

\medskip

\noi {\it Proof of Theorem \ref{thm1}.}  (i)
Follows from Lemmas \ref{lema2}  and \ref{lema4} (i) using standard argument.

\smallskip

\noi (ii) By Taylor's expansion,
\begin{equation*}
0 \ = \ \nabla L_n(\widehat \theta_n) = \nabla L_n(\theta_0) + \nabla^T \nabla L_n (\theta^*_n) (\widehat \theta_n - \theta_0),
\end{equation*}
where  $\theta^*_n \to_{\rm p} \theta_0  $ since   $\widehat \theta_n \to_{\rm p} \theta_0.  $
% $|\theta^*_n - \theta_0| \le |\widehat \theta_n - \theta_0| $.
Then $\nabla^T \nabla L_n (\theta^*_n)   \to_{\rm p} \nabla^T \nabla L (\theta_0)$
by Lemma \ref{lema4} \eqref{convL3}.
Next, since $\{r^2_t/\sigma^2_t(\theta_0) -1, {\cal F}_t, t \in \Z\} $ is a square-integrable and ergodic martingale difference sequence,
the convergence $n^{1/2}  \nabla L_n(\theta_0)   \stackrel{law}{\rightarrow} \ N(0, A(\theta_0)) $ follows
by the martingale central limit theorem in (\cite{Bill1968}, Thm.23.1). Then \eqref{clt1}
follows by Slutsky's theorem and \eqref{ABmat}. \hfill $\Box$

\medskip

\noi {\it Proof of Theorem \ref{thm2}.} Part (i)
follows from Lemmas \ref{lema2}  and \ref{lema4} (i) as in the case of Theorem \ref{thm1} (i).  To prove
part (ii), by Taylor's expansion
\begin{equation*}
0 \ = \ \nabla \widetilde L^{(\beta)}_n(\widetilde \theta^{(\beta)}_n) =
\nabla \widetilde L^{(\beta)}_n(\theta_0) + \nabla^T \nabla \widetilde L^{(\beta)}_n (\widetilde \theta^*_n) (\widetilde \theta^{(\beta)}_n - \theta_0),
\end{equation*}
where  $\widetilde \theta^*_n \to_{\rm p} \theta_0  $ since   $\widetilde \theta^{(\beta)}_n \to_{\rm p} \theta_0.  $  Then
$\nabla^T \nabla \widetilde L^{(\beta)}_n (\theta^*_n)   \to_{\rm p} \nabla^T \nabla L (\theta_0)$
by Lemma \ref{lema4} \eqref{convL3}- \eqref{convL4}. From the  proof of Theorem \ref{thm1} (ii) we have that
$n^{\beta/2}  \nabla L^{(\beta)}_n(\theta_0)   \stackrel{law}{\rightarrow} \ N(0, A(\theta_0)), $
where $L^{(\beta)}_n(\theta)  := \frac{1}{[n^\beta]}\sum_{t=n-[n^\beta]+1}^n  {l}_{t}(\theta). $ Hence,
the central limit theorem in \eqref{conv}  follows from
\begin{equation}\label{beta}
I_n(\beta) := \E |\nabla \widetilde L^{(\beta)}_n(\theta_0) - \nabla L^{(\beta)}_n(\theta_0)| \  = \ o(n^{-\beta/2}).
\end{equation}
We have $I_n(\beta) \le \sup_{n-[n^\beta] \le t \le n}  \E |\nabla {l}_{t}(\theta_0)-  \nabla \widetilde {l}_{t}(\theta_0)|$ and
\eqref{beta} follows from
\begin{equation}\label{beta1}
\E |\nabla {l}_{t}(\theta_0)-  \nabla \widetilde {l}_{t}(\theta_0)| \  = \ o(t^{-\beta/2}), \qquad t \to \infty.
\end{equation}
Write $\|\xi\|_p := \E^{1/p} |\xi|^p $ for $L^p$-norm of r.v. $\xi$.
Using $|\nabla ({l}_{t}(\theta_0)-  \widetilde {l}_{t}(\theta_0))| \le r^2_t |\nabla (\sigma^{-2}_{t}(\theta_0)-  \widetilde \sigma^{-2}_{t}(\theta_0))|
+ |\nabla ( \log \sigma^{2}_{t}(\theta_0)-  \log \widetilde \sigma^{2}_{t}(\theta_0))| $ and assumption $\E |r_t|^5 < \infty $,
relation
\eqref{beta1} follows from
\begin{eqnarray} \label{SS}
\|\sigma^{-4}_t \partial_i \sigma^2_t -  \widetilde \sigma^{-4}_t \partial_i \widetilde \sigma^2_t\|_{5/3}
&=&O(t^{d_0 -1/2} \log t) \quad \text{and} \\
\|\sigma^{-2}_t \partial_i \sigma^2_t -  \widetilde \sigma^{-2}_t \partial_i \widetilde \sigma^2_t\|_1
&=&O(t^{d_0-1/2}\log t), \qquad  i=1, \cdots, 5,  \nn
% \sigma^{-2}_t(\theta))|^{2} +
%\E^{1/2} |\nabla ( \log \sigma^{2}_{t}(\theta_0)-  \log \widetilde \sigma^{2}_{t}(\theta_0))|^2 \big)
\end{eqnarray}
where $\sigma^2_t := \sigma^2_t(\theta_0), \, \widetilde \sigma^2_t := \widetilde \sigma^2_t(\theta_0), \,
\partial_i \sigma^2_t := \partial_i \sigma^2_t(\theta_0), \, \partial_i \widetilde \sigma^2_t := \partial_i \widetilde \sigma^2_t(\theta_0)$.
Below, we prove the first relation \eqref{SS} only, the proof of the second one being similar.
We have  $\sigma^{-4}_t \partial_i \sigma^2_t -  \widetilde \sigma^{-4}_t \partial_i \widetilde \sigma^2_t
= \sigma^{-4}_t \widetilde \sigma^{-4}_t (\widetilde \sigma^2_t + \sigma^2_t) (\widetilde \sigma^2_t - \sigma^2_t) \partial_i \sigma^2_t
+ \widetilde \sigma^{-4}_t (\partial_i \sigma^2 - \partial_i \widetilde \sigma^2_t).$  Then using $\sigma^2_t \ge \omega^2_1/(1- \gamma_2) >0,  \widetilde \sigma^2_t \ge \omega^2_1/(1-\gamma_2) >0$, relation the first relation in \eqref{SS}  follows from
\begin{eqnarray} \label{SS1}
\|(\sigma^2_t - \widetilde \sigma^2_t)(\partial_i \sigma^2_t/\sigma_t)\|_{5/3}&=&O(t^{d_0-1/2}) \quad \text{and} \\
\|\partial_i \sigma^2_t -  \partial_i \widetilde \sigma^2_t\|_{5/3}&=&O(t^{d_0- 1/2}\log t), \qquad  i=1, \cdots, 5. \label{SS2}
\end{eqnarray}
Consider \eqref{SS1}.  By H\"older's inequality,  $\|(\sigma^2_t - \widetilde \sigma^2_t)(\partial_i \sigma^2_t/\sigma_t)\|_{5/3} \le
 \|\sigma^2_t - \widetilde \sigma^2_t\|_{5/2} \|\partial_i \sigma^2_t/\sigma_t\|_5 $,
where $\|\partial_i \sigma^2_t/\sigma_t\|_5 < C $ according to \eqref{Esup1}. Hence,   \eqref{SS1} follows from
\begin{eqnarray} \label{SS3}
\|\sigma^2_t - \widetilde \sigma^2_t\|_{5/2}&=&O(t^{d_0-1/2}).
\end{eqnarray}
To show \eqref{SS3}, similarly as in the proof of  \eqref{V11}
split $\sigma_t^2 -  \widetilde \sigma^2_t = U_{t,1}  + U_{t,2}$, where
$U_{t,i}:= U_{t,i}(\theta_0), i =1,2 $ are defined in \eqref{Udef}, i.e.,
$U_{t,1} = \sum_{\ell =1}^{t-1}  \gamma^\ell_0 \big\{\big(a_0 + c_0 Y_{t-\ell}\big)^2  - \big(a_0 + c_0 \widetilde Y_{t-\ell}\big)^2  \big\}, \,
U_{t,2} =\sum_{\ell =t}^{\infty}  \gamma^\ell_0 \big\{\omega^2_0 + \big(a_0 + c_0  Y_{t-\ell}\big)^2\big\} $ and
$Y_t := Y_t(d_0), \widetilde Y_t := \widetilde Y_t (d_0). $  We have
%\begin{eqnarray*}
$|U_{t,1}|\le C \sum_{\ell=1}^{t-1} \gamma_0^\ell |Y_{t-\ell} - \widetilde Y_{t-\ell}| (1 + |Y_{t-\ell}| +
|\widetilde Y_{t-\ell}|),   \
|U_{t,2}|\ \le\  C\sum_{\ell =t}^\infty \gamma^\ell_0 (1 + |Y_{t-\ell}|^2)$  %\end{eqnarray*}
and hence
\begin{eqnarray}
\|\sigma^2_t - \widetilde \sigma^2_t\|_{5/2}&\le&C\Big\{\sum_{\ell =1}^{t-1} \gamma_0^\ell \|(Y_{t-\ell} - \widetilde Y_{t-\ell}) (1 + |Y_{t-\ell}|+|\widetilde Y_{t-\ell}|)\|_{5/2}
+  \sum_{\ell=t}^\infty \gamma^\ell_0 (1 + \|Y_{t-\ell}\|_{5})\Big\} \nn \\
&\le&C\Big\{\sum_{\ell =1}^{t-1} \gamma_0^\ell \|Y_{t-\ell} - \widetilde Y_{t-\ell}\|_{5}
+  \sum_{\ell=t}^\infty \gamma^\ell_0 \Big\},   \label{SS4}
\end{eqnarray}
where we used the fact that $\|Y_t\|_5 < C $, $\|\widetilde Y_t\|_5 < C $ by $\|r_t\|_5 < C$  and
Rosenthal's inequality in \eqref{rosen}. In a similar way from  \eqref{rosen} it follows that
\begin{equation} \label{YY}
\|Y_{t-\ell} - \widetilde Y_{t-\ell}\|_{5} \ \le \ C \big\{\sum_{j > t-\ell} j^{2(d_0-1)}\big\}^{1/2} \ \le  \ C (t-\ell)^{d_0-1/2}.
\end{equation}
Substituting \eqref{YY}
into \eqref{SS4} we obtain
\begin{eqnarray*}
\|\sigma^2_t - \widetilde \sigma^2_t\|_{5/2}
&\le&C\Big\{\sum_{\ell =1}^{t-1} \gamma_0^\ell (t-\ell)^{d_0-1/2}
+  \sum_{\ell=t}^\infty \gamma^\ell_0 \Big\} \ = \ O(t^{d_0-1/2}),
\end{eqnarray*}
proving \eqref{SS3}.

It remains to show  \eqref{SS2}. %With $i=1$ the proof is analogous to \eqref{SS4} so we show only the case $i\neq 1$.
Similarly as above, $\partial_i \sigma^2_t -  \partial_i \widetilde \sigma^2_t = \partial_i U_{t,1} + \partial_i U_{t,2} $, where
$\partial_i U_{t,j} := \partial_i U_{t,j}(\theta_0), j=1,2$.
Then \eqref{SS2} follows from
\begin{equation} \label{Ui}
\|\partial_i U_{t,1}\|_{5/3}=O(t^{d_0-1/2}\log t) \quad \text {and} \quad
\|\partial_i U_{t,2}\|_{5/3}= o(t^{d_0-1/2}), \quad i=1, \cdots, 5.
\end{equation}
For $i=1$, the proof of \eqref{Ui} is similar to  \eqref{SS4}. Consider \eqref{Ui} for
$2\le i \le 5$.  Denote $V_{t}(\theta):=2a + c(Y_{t}(d) + \widetilde Y_{t}(d)), \, V_t := V_t(\theta_0), \partial_i V_t := \partial_i V_t(\theta_0)   $,
then
\begin{eqnarray*}
\|\partial_i U_{t,1}\|_{5/3}
%&\le&
%C\sum_{\ell =1}^{t-1} \gamma_0^\ell  \Big\{
%\|\partial_i \big((Y_{t-\ell}(\theta - \widetilde Y_{t-\ell}) V_t\big) \|_{5/3}
%+\|(Y_{t-\ell} - \widetilde Y_{t-\ell}) \partial_i (2a + Y_{t-\ell} + \widetilde Y_{t-\ell}) \|_{5/3}
%\Big\}\\
&\le&
C\sum_{\ell =1}^{t-1} \gamma_0^\ell  \Big\{
\|\partial_i (Y_{t-\ell} -  \widetilde Y_{t-\ell})\|_{5} \|V_t\|_5
+\|Y_{t-\ell} - \widetilde Y_{t-\ell} \|_{5}   \|\partial_i V_t\|_{5}
\Big\},
\end{eqnarray*}
where $\partial_i (Y_{t-\ell} -  \widetilde Y_{t-\ell}) = 0, \partial_i \ne \partial_d $ and
\begin{eqnarray*}
\|\partial_d (Y_{t} -  \widetilde Y_{t})\|_5
&=&\|\sum_{j> t} j^{d_0-1} (\log j) r_{t-j} \|_5 \\
&\le&C\big\{\sum_{j>t} j^{2(d_0-1)} \log^2 j \big\}^{1/2}  \ = \ O( t^{d_0 -1/2} \log t)
\end{eqnarray*}
similarly as in \eqref{YY}  above.
Hence, the first relation in \eqref{Ui}
%$\|\partial_i(U_{t,1})\|_{5/3}=O(t^{d_0-1/2})$
follows from \eqref{YY} and
$ \|\partial_i V_t \|_{5} \le C(1+ \|\partial_d Y_{t-\ell}  \|_{5}+  \| \partial_d \widetilde Y_{t-\ell} \|_{5} )\le C<\infty$
%and $\|\partial_i (Y_{t-\ell} - \widetilde Y_{t-\ell})\|_{5} =O(t^{d_0-1/2})$
as in the proof of \eqref{SS2}, and the proof of the second relation in \eqref{Ui} is analogous.
This proves \eqref{beta1} and
completes the proof of Theorem \ref{thm2}. \hfill $\Box$

\section*{Acknowledgement}

This research was supported in part by grant MIP-063/2013 from the Research Council of Lithuania.

\begin{table}
\centering
\caption {Sample RMSE of QML estimates of  $\theta_0 = (\gamma_0, \omega_0, a_0, c_0, d_0)$ of the GQARCH process in \eqref{sigma3}
%$\widehat \theta_n = ( \widehat \gamma_n, \widehat \omega_n, \widehat a_n, \widehat c_n, \widehat d_n)$,
for $\gamma_0 = 0.7, a_0=-0.2, c_0=0.2$ and different values of $\omega_0, d_0$. The number of
replications is $100$ %[\textbf{MATLAB iverciai}].
%based on $N=100$ replications
%with $\gamma_0 = 0.7, a_0=0.2, c_0=0.2$ for different values of $\omega_0$ and long memory parameter $d_0$.
}\label{tab:title}
\bigskip
\begin{tabular}{lp{3 em}ccccc}
\hline
& & \multicolumn{5}{c}{$\omega_0$=0.1} \\ \cline{3-7}
$n$ & $d_0$ &$\widehat \gamma_n$ & $\widehat \omega_n$ & $\widehat a_n$ & $\widehat d_n$ & \multicolumn{1}{c}{$\widehat c_n$} \\
\hline
1000& 0.1 & $0.091$ & $0.057$ & $0.035$ &  $0.103$  & $0.035$ \\
& 0.2 & $0.083$ & $0.047$ & $0.045$ &  $0.109$  & $0.031$ \\
& 0.3 & $0.071$ & $0.045$ & $0.047$ &  $0.094$  & $0.043$ \\
& 0.4 & $0.073$ & $0.029$ & $0.054$ &  $0.097$  & $0.036$ \\
\\
5000& 0.1 & $0.031$ & $0.021$ & $0.012$ &  $0.047$  & $0.015$ \\
& 0.2 & $0.030$ & $0.015$ & $0.015$ &  $0.041$  & $0.014$ \\
& 0.3 & $0.028$ & $0.011$ & $0.025$ &  $0.042$  & $0.013$ \\
& 0.4 & $0.031$ & $0.014$ & $0.053$ &  $0.059$  & $0.018$ \\
\hline
\end{tabular}
\\
\begin{tabular}{lp{3 em}ccccc}
\hline
& & \multicolumn{5}{c}{ $\omega_0$=0.01} \\ \cline{3-7}
$n$ & $d_0$ &$\widehat \gamma_n$ & $\widehat \omega_n$ & $\widehat a_n$ & $\widehat d_n$ & \multicolumn{1}{c}{$\widehat c_n$} \\
\hline
1000& 0.1 & $0.070$ & $0.049$ & $0.030$ &  $0.103$  & $0.029$ \\
& 0.2 & $0.061$ & $0.043$ & $0.035$ &  $0.089$  & $0.024$ \\
& 0.3 & $0.066$ & $0.040$ & $0.045$ &  $0.106$  & $0.044$ \\
& 0.4 & $0.055$ & $0.042$ & $0.056$ &  $0.105$  & $0.038$ \\
\\
5000& 0.1 & $0.025$ & $0.032$ & $0.011$ &  $0.035$  & $0.013$ \\
& 0.2 & $0.022$ & $0.028$ & $0.013$ &  $0.032$  & $0.013$ \\
& 0.3 & $0.025$ & $0.028$ & $0.025$ &  $0.046$  & $0.016$ \\
& 0.4 & $0.031$ & $0.031$ & $0.046$ &  $0.096$  & $0.034$ \\
\hline
\end{tabular}
\\

\begin{tabular}{lp{3 em}ccccc}
\hline
& & \multicolumn{5}{c}{ $\omega_0$=0.001} \\ \cline{3-7}
$n$ & $d_0$ &$\widehat \gamma_n$ & $\widehat \omega_n$ & $\widehat a_n$ & $\widehat d_n$ & \multicolumn{1}{c}{$\widehat c_n$} \\
\hline
1000& 0.1 & $0.086$ & $0.058$ & $0.026$ &  $0.095$  & $0.037$ \\
& 0.2 & $0.056$ & $0.043$ & $0.027$ &  $0.084$  & $0.031$ \\
& 0.3 & $0.053$ & $0.039$ & $0.046$ &  $0.080$  & $0.029$ \\
& 0.4 & $0.055$ & $0.047$ & $0.060$ &  $0.122$  & $0.041$ \\
\\
5000& 0.1 & $0.022$ & $0.033$ & $0.009$ &  $0.031$  & $0.012$ \\
& 0.2 & $0.020$ & $0.030$ & $0.012$ &  $0.028$  & $0.012$ \\
& 0.3 & $0.022$ & $0.032$ & $0.024$ &  $0.038$  & $0.014$ \\
& 0.4 & $0.032$ & $0.037$ & $0.046$ &  $0.098$  & $0.031$ \\
\hline
\end{tabular}
\end{table}

\end{document}